\documentclass[10pt]{amsart}

\input{SROP.sty}

\title[Realization of some Stanley-Reisner algebras and graph colorings]{Realization of some Stanley-Reisner algebras \\ and graph colorings}
\author{Yang Hu and Donald Stanley}
\address{Department of Mathematics and Statistics, University of Regina}
\email{Yang.Hu@uregina.ca}
\address{Department of Mathematics and Statistics, University of Regina}
\email{Donald.Stanley@uregina.ca} 
\date{}         


\begin{document}

\begin{abstract}
It is a classical problem in algebraic topology to decide whether a given graded $\ZZ$-algebra can be realized as the integral cohomology ring of a space. In this paper, we introduce families of Stanley-Reisner algebras depending on graphs and relate their realizability to the span coloring and chromatic number of the graph.
\end{abstract}

\maketitle

\tableofcontents

\section{Introduction} \label{sec:intro}

The cohomology realization problem \cite{Ste61}, asks whether a given graded $\ZZ$-algebra can be the cohomology ring of a space. One notable example is the realizability of the truncated polynomial algebra $\ZZ[x]/(x^3)$. This is essentially the Hopf invariant one problem, which was resolved in the ground-breaking work of Adams \cite{Ada60} using cohomology operations. Such an algebra is proved to be realizable precisely when the degree of the generator equals $2$, $4$, or $8$. 
Although the realization problem seems intractable in its full generality, certain special cases have been  extensively studied. These include the case of polynomial algebras, whose realization is completely determined in  \cite[Theorem 1.1]{AG08}. The realizable ones are precisely those isomorphic to a tensor product of 
\begin{itemize}
 \item[-] $\Hmr^*(\CP^{\infty}; \ZZ) \cong \ZZ[x]$, $|x| = 2$,
 \item[-] $\Hmr^*(\BSU(n); \ZZ) \cong \ZZ[x_2, \cdots, x_n]$, $|x_i| = 2i$ for $i = 2, \cdots, n$, and
 \item[-] $\Hmr^*(\BSp(m); \ZZ) \cong \ZZ[y_1, \cdots, y_m]$, $|y_j| = 4j$ for $j = 1, \cdots, m$.
\end{itemize}

\vspace{1mm}

Recent advancements in the resolution of the cohomology realization problem feature the realization of certain monomial ideal rings that are combinatorial in flavor. Recall that a simplicial complex $\K$ on a finite vertex set $\V$ is a subset of the power set $\Pcal(\V)$ which is closed under taking subsets. Simplicial complexes are natural indexing objects for many important combinatorial constructions in topology, such as polyhedral products, and their special cases including moment-angle complexes and Davis-Januszkiewicz spaces (see \cite{BBC20} for a survey of these constructions).
On the algebraic side, given a simplicial complex $\K$ together with a positive even grading $\upphi: \V\rightarrow 2\ZZ_{>0}$ on its vertex set $\V = \{v_1, \cdots, v_n\}$, there is an associated monomial ideal ring called the Stanley-Reisner ring (also known as the face ring):
\[
\SR(\K, \upphi):=\ZZ[v_1, \cdots, v_n]/I_{\K},
\]
where the monomial ideal $I_{\K}$ is generated by monomials corresponding to simplicies that do not appear in $\K$. In the context of Steenrod's problem, families of Stanley-Reisner rings are known to admit topological realizations. 
For example, when every vertex of $\K$ has degree $2$, the Stanley-Reisner ring $\SR(\K, \upphi)$ can be realized by the Davis-Januszkiewicz space $DJ(\K)$ \cite{DJ91}. When the vertices of $\K$ are of degree $2$ or $4$, the work of Bahri, Bendersky, Cohen and Gitler establishes the realizability of $\SR(\K, \upphi)$ by polyhedral products \cite[Theorem 2.35]{BBCG10}. Recent work of Takeda \cite{Tak24} provides necessary and sufficient conditions for the realizability of fairly general Stanley-Reisner rings, where the $\SR(\K, \upphi)$ in question must satisfy the following property: for any $x, y \in \SR(\K, \upphi)$ whose degrees are powers of $2$, the product $xy$ must be zero in $\SR(\K, \upphi)$.

\vspace{1mm}

In a recent preprint of Stanley and Takeda \cite{ST25}, the authors study graph-related Stanley-Reisner rings at the other end of the spectrum, where the generators are in degree $4$ or $6$ and there are no relations between the degree $4$ generators. More precisely, given a finite graph $\G$ (regarded as a one-dimensional simplicial complex), the authors define
\begin{equation}\label{informal_algebra}
\A(n, \G):=\SR(\K, \upphi) = \ZZ[x_1, \cdots, x_n, y_1, \cdots, y_m]/I_{\K}
\end{equation}
where $\K = \sfdel^{n-1}*\G$ is the join of the standard $(n-1)$-simplex $\sfdel^{n-1}$ with $\G$ (recall that the simplicial join respects topological realization). The degree assignment $\upphi$ is such that $\upphi(x_i) = 4$ for vertices $x_i$ from $\sfdel^{n-1}$, and $\upphi(y_j)=6$ for vertices $y_j$ from the graph $\G$.
If $\A(n, \G)$ is realizable, then $\A(n, \G)\otimes \ZZ/p$ necessarily supports an action of the mod-$p$ Steenrod algebra $\Acal_p$, for every prime $p$. In \cite{ST25}, the authors develop an algebraic version of graph coloring -- called the span coloring -- and use it to detect when $\A(n, \G)\otimes \ZZ/2$ has an action of $\Acal_2$. 
Recall that a usual $n$-coloring of a graph $\G$ is a function $f$ from the vertex set $\V(\G)$ to an $n$-element set, say $\{1, \cdots, n\}$, so that $f(v_i)\neq f(v_j)$ whenever two vertices $v_i$ and $v_j$ are neighboring in $\G$ (i.e., connected by an edge). The chromatic number of $\G$, written as $\chi(\G)$, is the smallest $n$ so that $\G$ admits an $n$-coloring. 
Given a field $\mathbf{k}$, an $n$-span coloring of $\G$ can be considered as a function $f: \V(\G) \rightarrow \mathbf{k}^n \setminus \{0\}$ so that $f(v_i)\notin \left\langle f(N(v_i)) \right\rangle$ for every $v_i\in \V(\G)$. Here $N(v_i)$ is the set of neighbors of $v_i$ and $\left\langle f(N(v_i)) \right\rangle$ denotes the sub $\mathbf{k}$-vector space spanned by vectors $f(v_j)$ for $v_j\in N(v_i)$. The span chromatic number of $\G$, denoted by $s_{\mathbf{k}}\chi(\G)$, is the smallest $n$ so that $\G$ admits an $n$-span coloring. When $\mathbf{k}$ is the finite field $\FF_p$, we write $s_{\mathbf{k}}\chi(\G)$ simply as $s_p\chi(\G)$. 

\vspace{1mm}

A graph-theoretic characterization of the admission of an $\Acal_2$-algebra structure on $\A(n, \G)\otimes \ZZ/2$ follows.

\begin{thm}[\cite{ST25}, Theorem 7.8] \label[thm]{ST1}
The algebra $\A(n, \G)\otimes \ZZ/2$ has an $\Acal_2$-action precisely when $s_2\chi(\G)\leq n$.
\end{thm}

In particular, if $\A(n, \G)$ is realizable then $s_2\chi(\G)\leq n$. The authors introduce a topologically defined graph invariant, $\chi_{\Top}(\G)$, as the smallest $n$ so that $\A(n, \G)$ is realizable. \Cref{ST1} can then be interpreted as proving $s_2\chi(\G)$ as a lower bound of $\chi_{\Top}(\G)$. It turns out that the usual chromatic number $\chi(\G)$ serves as an upper bound:

\begin{thm}[\cite{ST25}, Theorem 8.4] \label[thm]{ST2}
If $\chi(\G) \leq n$, then the algebra $\A(n, \G)$ is realizable.
\end{thm}

Combining \Cref{ST1} and \Cref{ST2}, one obtains:

\begin{thm}[\cite{ST25}, Theorem 8.6] \label[thm]{ST3}
For any graph $\G$, $s_2\chi{\G} \leq \chi_{\Top}(\G) \leq \chi(\G)$.
\end{thm}

\vspace{1mm}

This paper is aimed at generalizing the above main results of \cite{ST25} to a wider range of Stanley-Reisner rings, and exhibiting connections at other primes. In \Cref{sec:prime3} we introduce the Stanley-Reisner algebra $\B(n, \G)$ (see \Cref{BnG}), which is the $3$-primary analogue of the algebra $\A(n, \G)$, with generators now in degrees $4$ and $8$. Our first main result is:

\begin{main}[See \Cref{thm:main3}] \label[main]{main1}
If there is an $\Acal_3$-algebra structure on $\B(n, \G) \otimes \ZZ/3$, then $s_3\chi(\G) \leq n$.
\end{main}

Note that \Cref{main1} only partly generalizes \Cref{ST1}. Indeed, at the prime $3$, having $s_3\chi(\G) \leq n$ need not guarantee the existence of an action of $\Acal_3$ on $\B(n, \G) \otimes \ZZ/3$, with the possible existence of further obstructions. 

When we consider realizations of $\B(n,\G)$ with extra structure we get that a solution to the realization problem is equivalent to $\chi(\G)\leq n$. 

\begin{main} \label{thm:lastthmintro}
The following are equivalent:
\begin{enumerate}[label=(\alph*)]
\item The chromatic number $\chi(\G)$ is less than or equal to $n$.
\item There is a realization $\Xmr_{\B(n, \G)}$ of $\B(n, \G)$ together with a map
\[
\prod_{i=1}^n \BSp(1)^{\times n} \longrightarrow \Xmr_{\B(n, \G)}
\]
realizing the standard projection
\[
p: \B(n, \G)\otimes \ZZ/3 \longrightarrow \ZZ/3[x_1, \cdots, x_n] \cong \Hmr^*(\BSp(1)^{\times n}).
\]
\item There is an $\Acal_3$-algebra structure on $\B(n, \G)\otimes \ZZ/3$ so that the projection $p$ is a map of  $\Acal_3$-algebras.
\end{enumerate}
\end{main}

\vspace{1mm}



In \Cref{sec:primep} we consider a general family $\A_p(\vb{s}, \G)$ of Stanley-Reisner algebras, depending on an odd prime $p$ and sequence of non-negative integers $\vb{s}=(s_1, \cdots , s_{p-1})$. These algebras are modeled on products of classifying spaces $\BSU(k)$ and $\BSp(l)$, generalizing both $\A(n, \G)$ at the prime $2$ and the $\B(n, \G)$ at the prime $3$ and $s_i$ gives the number of generators of degree $2i+2$. (See \Cref{ApsG}.) We prove:

\begin{main}[See \Cref{thm:mainAp}] \label[main]{main3}
If there is an $\Acal_p$-algebra structure on $\A_p(\vb{s}, \G) \otimes \ZZ/p$, then $s_p\chi(\G) \leq s_1$, where $s_1$ is the first coordinate of the vector $\vb{s}$.
\end{main}

For the next theorem we will 
assume that $\vb{s}=\vb{r}+\vb{t}$ where 
$\vb{r}=(r_1,\cdots,r_{p-1})$ and $\vb{t}=(t_1,\cdots,t_{p-1})$ are sequences such that $r_i\geq r_{i+1}$, 
$t_{2i}=0$ and $t_i\geq t_{i+2}$. The first sequence corresponds to products of various $BSU(n)$ and the second to products of $BSp(n)$. 

\begin{main}[See \Cref{thm:realeyes}] \label[main]{main4}
If $r_{p-1}+t_{p-2}\leq \chi(\G)$ then the algebra $\A_p(\vb{s}, \G)$ is realizable. 
\end{main}

Given an integer $n$ and taking $\vb{r}(n)$ as the constant sequence $n$, we get 
a topologically defined graph invariant using the algebra $\A_p(\vb{r}(n), \G)$. Namely, we define $\chi_{\Top_p, \A}(\G)$ as the smallest $n$ so that $\A_p(\vb{r}(n), \G)$ is realizable. The next result follows immediately, generalizing \Cref{ST3}.


\begin{main}[See \Cref{thm:chronumber}] \label[main]{main5}
For any graph $\G$, we have $s_p\chi(\G) \leq \chi_{\Top, \A}(\G) \leq \chi(\G)$.
\end{main}

\subsection*{Acknowledgment}
This work is supported by the Natural Sciences and Engineering Research Council of
Canada Grant RGPIN-05466-2020. The first author would like to thank the Pacific Institute for the Mathematical Sciences (PIMS) for their support when conducting research for this paper. PIMS Report Identifiers: PIMS-20251030-CRG41; PIMS-20251030-PDF.

\vspace{2mm}

\section{The Stanley-Reisner algebra \texorpdfstring{$\B(n, \G)$}{BnG}} \label{sec:prime3}

\begin{defn} \label[defn]{BnG}
Let $\G$ be a finite graph. Write $\B(n, \G)$ for the Stanley-Reisner ring $\SR(\K, \upphi)$, where $\K$ stands for the simplicial complex $\sfdel^{n-1} * \G$, so that
\[
\SR(\K, \upphi) = \ZZ[x_1, \cdots, x_n, y_1, \cdots, y_m]/I_{\K}, \quad \upphi(x_i) = 4 \text{ and } \upphi(y_j) = 8.
\]
\end{defn}

Our goal is to prove the following result.

\begin{thm} \label{thm:main3}
If there is an $\Acal_3$-algebra structure on $\B(n, \G) \otimes \ZZ/3$, then $s_3\chi(\G) \leq n$.
\end{thm}

We need the following two technical lemmas from \cite{ST25}.

\begin{lemma}[\cite{ST25}, Lemma 6.1] \label[lemma]{preserve}
Let $p$ be a prime, and assume that $\SR(\K, \upphi)\otimes \ZZ_p$ is already an $\Acal_p$-algebra. Then for $\sigma \in P_{\max}(\K)$, the ideal $(\Vmr \setminus \sigma)$ is preserved by the $\Acal_p$-action.
More generally, if $\Umr$ is a subset of $P_{\max}(\K)$ then the ideal $\cap_{\tau\in \Umr} (\Vmr \setminus \tau)$ is preserved by the $\Acal_p$-action.
\end{lemma}

\begin{lemma}[\cite{ST25}, Lemma 6.2] \label[lemma]{ideal}
Let $\Umr \subset \K$ be the collection of all maximal simplicies. Then $\K = \cup_{\sigma\in \Umr} \sigma$, and $I_{\K} = \cap_{\sigma \in \Umr} (\Vmr \setminus \sigma)$. Suppose $\Umr'\subset \Umr$, and $\K':=\cup_{\sigma\in \Umr'}\sigma \subset \K$. Then $I = \cap_{\sigma\in \Umr'}(\Vmr\setminus \sigma)$ is the kernel of $\SR(\K, \upphi) \rightarrow \SR(\K', \upphi)$.
\end{lemma}

Assume that $\B(n, \G)\otimes \ZZ_3$ is an $\Acal_3$-algebra. We are interested in learning the behavior of $\Sp^1$ and $\Sp^3$ on $y_i$. 

\begin{prop} \label[prop]{square_formulas}
Suppose that $\B(n, \G)\otimes \ZZ_3$ has an action of $\Acal_3$, and that $\G$ is a graph where every vertex has degree at least two. Then:
\begin{enumerate}[label=(\arabic*)]
 \item The element $\Sp^1(y_i)$ is divisible by $y_i$.
 \item The element $\Sp^3(y_i)$ can be expressed as
\[
y_i^2g(y_i) + y_ih(y_i) + \sum_{1\leq j < k \leq m}\mu^{(i)}_{j, k}y_jy_k, \quad \text{where}
\] 
\begin{itemize}
\item[-] $g(y_i)$ is a linear combination of $x_1, \cdots, x_n$.
\item[-] $h(y_i)$ is a linear combination of length-three monomials $x_ax_bx_c$, $a, b, c\in \{1, \cdots, n\}$.
\item[-] $\mu^{(i)}_{j, k}$ is a linear combination of $x_1, \cdots, x_n$.
\end{itemize}
\end{enumerate}
\end{prop}

\begin{proof}
Our argument is similar to that in \cite[Lemma 6.7]{ST25}. Let $\Umr$ be the collection of maximal simplicies not containing $y_i$. Since each vertex of $\G$ has degree at least 2, for any $j\neq i$ there is some simplex $\tau \in \Umr$ so that $\tau$ contains the vertex $j$ and does not contain $i$. This implies that $\cup_{\tau\in \Umr} \tau = \K_{\V(\Kmr)\setminus \{i\}}$, and it follows from \Cref{ideal} that
\[
\cap_{\sigma \in \Umr} (\Vmr \setminus \sigma) = I_{\cup_{\tau\in \Umr} \tau} \subset (y_i) + (y_jy_k: 1\leq j < k \leq m).
\]
It then follows from \Cref{preserve} that  $\Sp^1(y_i)$ lies in $(y_i) + (y_jy_k: 1\leq j < k \leq m)$.
Since $\Sp^1(y_i)$ is of degree $12$, it must be of the form $y_if(y_i)$ where $f(y_i)$ is a $\ZZ/3$-linear combination of $x_1, \cdots, x_n$.
The same argument also implies that $\Sp^3(y_i)\in (y_i) + (y_jy_k: 1\leq j < k \leq m)$. Now that $\Sp^3(y_i)$ is of degree $20$, it must be of the form described above.
\end{proof}

Note that the values of $g(y_i)$ can be packaged to define a function
\[
g: \{y_1, \cdots, y_m\} \longrightarrow \ZZ/3\langle x_1, \cdots, x_n \rangle, \quad y_i \longmapsto g(y_i).
\]

We show that $g$ gives rise to a $n$-span coloring of $\G$.

\begin{lemma} \label[lemma]{P3_2G_coloring}
Assume that $\B(n, \G)\otimes \ZZ_3$ is an $\Acal_3$-algebra, and that every vertex of $\G$ has degree at least $2$. The linear inclusion
\[
\langle g(y_j) : y_j\in N(y_i) \rangle \longrightarrow \langle g(y_j) : y_j\in N(y_i) \cup \{y_i\} \rangle
\]
has nonzero cokernel for every $i$.
\end{lemma}

\begin{proof}
Assume for contradiction that it does not. Then $g(y_i)$ is a linear combination 
\begin{equation} \label{eq:assumption}
g(y_i) = \sum_{y_j\in N(y_i)} a_jg(y_j), \,\, a_j\in \ZZ/3.
\end{equation}
We now apply the Adem relation $\Sp^1\Sp^3 = \Sp^4$ to $y_i$, and use the knowledge of $\Sp^1(y_i)$ and $\Sp^3(y_i)$ from \Cref{square_formulas}. One obtains
\begin{equation} \label{eq:AdemP4}
y_i^3 = \Sp^4(y_i) = \Sp^1(\Sp^3(y_i)) = \Sp^1\left(y_i^2g(y_i) + y_ih(y_i) + \sum_{1\leq j < k \leq m}\mu^{(i)}_{j, k}y_jy_k\right).
\end{equation}
The right-hand side of \Cref{eq:AdemP4} must contain $y_i^3$ as a summand, dictating that $\Sp^1(g(y_i))$ contains $y_i$ as a summand. Now applying $\Sp^1$ to \Cref{eq:assumption}, one concludes that there exists some $j$ with $y_j\in N(y_i)$, so that $\Sp^1(g(y_j))$ contains $y_i$ as a summand. For such $j$, we find that
\begin{equation} \label{eq:AdemP4j}
y_j^3 = \Sp^4(y_j) = \Sp^1\Sp^3(y_j) = \Sp^1\left( y_j^2g(y_j) + y_jh(y_j) + \sum_{1\leq k < l \leq m}\mu^{(j)}_{k, l}y_ky_l \right)
\end{equation}
Since $\Sp^1(g(y_j))$ contains $y_i$ as a summand, the right-hand side of \Cref{eq:AdemP4j} contains $y_j^2y_i$ as a summand. Note that $y_j^2y_i$ is nonzero as $y_j\in N(y_i)$. However, $y_j^2y_i$ does not appear in the left-hand side of \Cref{eq:AdemP4j}, which is a contradiction.
\end{proof}

We are now ready to prove the main result, \Cref{thm:main3}.

\begin{proof}[Proof of \Cref{thm:main3}]
Let $2\G$ be the maximal subgraph of $\G$ whose vertices all have degree at least two.
There is a finite sequence $\G = \G_0 \supset \G_1 \supset \cdots \supset \G_{N-1} \supset \G_N = 2\G$ of subgraphs where each step removes a vertex of degree 0 or 1. (See \cite[Lemma 6.10]{ST25} for a proof.)
By \cite[Lemma 6.11]{ST25}, each step $\B(n, \G_{i+1})\otimes \ZZ/3 \longrightarrow \B(n, \G_i)\otimes \ZZ/3$ induces an $\Acal_3$-algebra structure on the target. 
It follows that $\B(n, \G_N)\otimes \ZZ/3 = \B(n, 2\G) \otimes \ZZ/3$ has an $\Acal_3$-action, and consequently that $s_3\chi(2\G) \leq n$ by \Cref{P3_2G_coloring}. 
It remains to show that $s_3\chi(\G) \leq n$, and we discuss two cases.

\vspace{1mm}

Suppose that $2\G$ is non-empty. Recall that there is a graph $\sfA_{(\ZZ/3)^n}$ representing $n$-span colorings over $\ZZ/3$, in the sense of \cite[Proposition 4.2]{ST25}.
Now that $s_3\chi(2\G) \leq n$, there is a graph map $\lambda: 2\G\longrightarrow \sfA_{(\ZZ/3)^n}$. 
Consider the finite sequence of graph inclusions: $2\G=\G_N\subset \G_{N-1} \subset \cdots \subset \G_1 \subset \G_0=\G$. Any graph map from $\G_{i+1}$ to $\sfA_{(\ZZ/3)^n}$ can be extended over $\G_i$ in the same way as explained in \cite[Theorem 6.12]{ST25}. 
So eventually, $\lambda$ can be extended to a map $\tilde{\lambda}: \G \longrightarrow \sfA_{(\ZZ/3)^n}$, which corresponds to an $n$-span coloring by \cite[Proposition 4.2]{ST25}. So $s_3\chi(\G) \leq n$ holds in this case.

\vspace{1mm}

Now suppose that $2\G$ is empty. In this case $\G$ is a union of trees, and therefore $s_3\chi(\G)\leq 2$. So if $n\geq 2$ then automatically $s_3\chi(\G) \leq 2 \leq n$. It then suffices to prove the result for $n=0$ and $n=1$. 
When $n=0$, we need to show that $s_3\chi(\G) = 0$, i.e., $\G$ is empty. Assume for contradiction that $\G \neq \emptyset$. Then a maximal simplex $\sigma$ in $\K = \sfdel^{n-1}*\G$ is either a vertex or an edge of $\G$. By \cite[Proposition 6.4]{ST25}, the algebra $\ZZ/3[\sigma]$ has an induced $\Acal_3$-action. That is to say, 
\[
A:= \ZZ/3[y], \, |y| = 8 \quad \text{or} \quad A':= \ZZ/3[y_1, y_2], \, |y_1| = |y_2| = 8
\]
has to support an $\Acal_3$-action. We show that this is impossible.
Indeed, if $A$ was an $\Acal_3$-algebra, then the Adem relation $\Sp^4 = \Sp^1\Sp^3$ gives the contradiction that
\[
y^3 = \Sp^4(y) = \Sp^1\Sp^3(y) = 0.
\]
(Here $\Sp^1\Sp^3(y) = 0$ simply for degree reasons.) Similarly, applying the same Adem relation to $y_1\in A'$ shows that $A'$ cannot admit an $\Acal_3$-action. 
When $n=1$, we need to show that $s_3\chi(\G) \leq 1$, i.e., $\G$ contains no edges.
Assume not. Then a maximal simplex $\sigma$ in $\K = \sfdel^{n-1}*\G$ is the join of an edge of $\G$ with the unique point in $\sfdel^0$. Again by \cite[Proposition 6.4]{ST25}, there is an induced $\Acal_3$-action on $\ZZ/3[\sigma]$, which is the algebra
\[
B:= \ZZ/3[x, y_1, y_2], \,\, |x|=4, |y_1| = |y_2| = 8.
\]
We find that $y_1^3 = \Sp^4(y_1) = \Sp^1\Sp^3(y_1)$, where $\Sp^3(y_1)$ can be a linear combination
\[
ax^5 + b_1x^3y_1 + b_2x^3y_2 + c_1xy_1^2 + c_2y_2^2 + c_3xy_1y_2.
\]
Applying the Cartan formula, one obtains
\[
\Sp^1\Sp^3(y_1) \equiv c_1\Sp^1(x)y_1^2 + c_1\Sp^1(x)y_2^2 + c_3\Sp^1(x)y_1y_2 \mod (x^2).
\]
Since $y_1^3$ must appear as a summand in the right-hand side of the above equation, the formula of $\Sp^1(x)$ must contain $y_1$ as a summand. By symmetry, $\Sp^1(x)$ also contains $y_2$ as a summand.
In other words,
\[
\Sp^1(x) = \lambda  x^2 + y_1 + y_2 \text{ for some } \lambda \in \ZZ/3.
\]
It follows that
\begin{align*}
y_1^3 &= \Sp^1\Sp^3(y_1) \\
 &\equiv c_1( \lambda  x^2 + y_1 + y_2)y_1^2 + c_1( \lambda  x^2 + y_1 + y_2)y_2^2 + c_3( \lambda  x^2 + y_1 + y_2)y_1y_2 \mod (x^2) \\
 &\equiv y_1^3 + y_2^3 + 2y_1^2y_2 + 2y_1y_2^2 \mod (x^2),
\end{align*}
which implies that $y_2^3 + 2y_1^2y_2 + 2y_1y_2^2$ is divisible by $x^2$. A contradiction.
\end{proof}

As a direct consequence of \Cref{thm:main3}, one obtains:

\begin{cor}
If $\B(n, \G)$ is realizable, then $s_3\chi(\G) \leq n$.
\end{cor}

\section{Characterizing the realizability of \texorpdfstring{$\B(n, \G)$}{BnG} }\label{subsec:discuss}

In this section we present more precise realization results for $\B(n, \G)$. Namely, $\B(n, \G)$ can be realized with some extra structure if and only if $\chi(\G)\leq n$ (see \Cref{thm:lastthm}). 
Recall that $\B(n, \G)$ is the Stanley-Reisner algebra
\[
\SR(\K, \upphi) = \ZZ[x_1, \cdots, x_n, y_1, \cdots, y_m]/I_{\K},
\]
where $\K = \sfdel^{n-1} * \G$, $m = |\G|$, $\upphi(x_i) = 4$  and $\upphi(y_j) = 8$. 


\begin{thm} \label{thm:lastthm}
The followings are equivalent:
\begin{enumerate}[label=(\alph*)]
\item The chromatic number $\chi(\G)$ is less than or equal to $n$.
\item There is a realization $\Xmr_{\B(n, \G)}$ of $\B(n, \G)$ together with a map
\[
\prod_{i=1}^n \BSp(1)^{\times n} \longrightarrow \Xmr_{\B(n, \G)}
\]
realizing the standard projection
\[
p: \B(n, \G)\otimes \ZZ/3 \longrightarrow \ZZ/3[x_1, \cdots, x_n] \cong \Hmr^*(\BSp(1)^{\times n}).
\]
\item There is an $\Acal_3$-algebra structure on $\B(n, \G)\otimes \ZZ/3$ so that the projection $p$ is a map of  $\Acal_3$-algebras.
\end{enumerate}
\end{thm}

\begin{rmk}
There is a similar result in the $2$-primary case for the Stanley-Reisner ring $\A(n, \G)$, already addressed in \cite{ST25}. In fact, for $\A(n, \G)$ the implications (a) $\Rightarrow$ (b) $\Rightarrow$ (c) always hold. When the number of vertices in $\G$ is at least three, we note that ${\sfdel^{n-1}} \in P_{\max} \K$, and the condition (c) is then equivalent to $s_2\chi(\G)\leq n$ by \cite[Theorem 7.8]{ST25}. The implication (c) $\Rightarrow$ (a) does not hold in general. A counter-example occurs when $\G$ is the graph $\A_{(\ZZ/2)^{3}}$, whose $s_2\chi(\G)$ is strictly smaller than $\chi(\G)$ (see \cite[Example 6.5]{ST25}).
\end{rmk}

The following technical lemma is key to the proof of \Cref{thm:lastthm}.

\begin{lemma} \label[lemma]{last}
Suppose that $2\G = \G$, and that $\B(n, \G)\otimes \ZZ/3$ admits an $\Acal_3$-action so that
\begin{equation} \label{eq:assume}
\Sp^1(x_i) = x_i^2 + g(x_i)
\end{equation}
for every $i=1, \cdots, n$, where $g(x_i)$ is some $\ZZ/3$-linear combination of the $y_j$'s.
Then the graph $\G$ can be $n$-colored, and therefore $\chi(\G)\leq n$.
\end{lemma}

\begin{proof} Since $2\G  = \G$, it follows from \Cref{square_formulas} that
\begin{equation} \label{eq:e1}
 \Sp^1(y_i) = y_ih(y_i)
\end{equation}
and that
\begin{equation} \label{eq:e2}
\Sp^3(y_i) = y_i^2f(y_i) + y_i\lambda(y_i) + \sum_{1\leq j<k\leq m} \mu^{(i)}_{j, k}y_jy_k 
\end{equation}
where $f(y_i), h(y_i), \mu^{(i)}_{j, k}$ are linear combinations of the $x_l$'s, and $\lambda(y_i)$ is a linear combination of length-three $x$-monomials $x_ax_bx_c$.
Also note that
\begin{equation} \label{eq:e3}
 \Sp^1\Sp^1(x_i) = 2\Sp^2x_i = 2x_i^3
\end{equation}
and that
\begin{equation} \label{eq:e4}
  \Sp^1\Sp^3(y_i) = \Sp^4y_i = y_i^3
\end{equation}

Fix an arbitrary $i\in \{1, \cdots, m\}$. Assuming $2\G = \G$, the algebra 
\[
\A_i:=\ZZ/3[x_1, \cdots, x_n, y_i]
\]
inherits an $\Acal_3$-action, with the projection
\[
\B(n, \G)\otimes \ZZ/3 \longrightarrow A_i
\]
being a map of $\Acal_3$-algebras. Therefore, \eqref{eq:assume}, \eqref{eq:e1}, \eqref{eq:e2}, \eqref{eq:e3}, and \eqref{eq:e4} still hold in $\A_i$. 
Consider the composite
\[
g_i: \{x_1, \cdots, x_n\} \stackrel{g}{\longrightarrow} \ZZ/3\langle y_1, \cdots y_m \rangle \longrightarrow  \ZZ/3\langle y_i \rangle.
\]
We claim that $g_i$ is nonzero (as a function). Indeed, by \eqref{eq:e2} and \eqref{eq:e4}
\begin{align*}
y_i^3 = \Sp^1\Sp^3y_i &= \Sp^1\left(y_i^2f(y_i) + y_i\lambda(y_i) + \sum_{1\leq j<k\leq m} \mu^{(i)}_{j, k}y_jy_k\right) \\
&= 2y_i\Sp^1(y_i)f(y_i)+y_i^2\Sp^1\left(f(y_i)\right) + \cdots.
\end{align*}
Note that the term $y_i^3$ can only possibly occur from the right-hand side by appearing as a summand in $V$. More precisely, suppose that $f(y_i) = c_1x_1+\cdots+c_nx_n$. Then
\begin{align*}
y_i^2\Sp^1\left(f(y_i)\right) &= y_i^2(c_1\Sp^1x_1 + \cdots + c_n\Sp^nx_n) \\
&= y_i^2(c_1x_1^2 + c_1g(x_1) + \cdots + c_nx_n^2 +c_ng(x_n)).
\end{align*}
So there exist some $j$ so that $g_i(x_j)$ is nonzero. Without loss of generality, assume
\[
g_i(x_1) = 2y_i.
\]

Now working in $\A_i$, \eqref{eq:e3} implies that
\[
2x_1^3 = \Sp^1\Sp^1x_1 = \Sp^1(x_1^2+2y_i) = 2x_1\Sp^1(x_i)+2\Sp^1(y_i) = 2x_1(x_1^2+2y_i)+2y_ih(y_i),
\]
and hence that
\[
h(y_i) = x_1
\]
Since the ``$i$'' we fix is arbitrary,  the map
\[
h: \{y_1, \cdots, y_m\} \longrightarrow \ZZ/3\langle x_1, \cdots, x_n \rangle
\]
factors through
\[
\hbar: \{y_1, \cdots, y_m\} \longrightarrow \{x_1, \cdots, x_n\}.
\]

We prove \Cref{last} by showing that $\hbar$ is a coloring. 
Working in $\A_i$ with $i\neq 1$, we have $g(y_i) = c_iy_i$ for some $c_i \in \ZZ/3$.
The Adem relation $\Sp^1\Sp^1 = 2\Sp^2$ implies that $c_i = 0$.
Furthermore, the Adem relation $\Sp^1\Sp^3 = \Sp^4$ implies that
\begin{align*}
 y_i^3 = \Sp^1\Sp^3y_i &= \Sp^1\left( y_i^2f(y_i) + y_i\lambda(y_i) \right) \\
  &= 2y_i\Sp^1(y_i)f(y_i) + y_i^2\Sp^1(f(y_i)) + \cdots \\
  &= 2y_i^2x_1f(y_i) + y_i^2(c_1x_1^2 + c_1g(x_1) + \cdots + c_nx_n^2 + c_ng(x_n)) + \cdots,
\end{align*}
where $f(y_i) = c_1x_1 + \cdots + c_nx_n$. Recall that $g(x_1) = 2y_i$, so $c_1 = 2$ and $c_j=0$ for $j>1$. In other words,
\[
f(x_i)  =2x_1.
\]
Now suppose that $j, k$ are connected by an edge in $\G$. If $\hbar(y_j) \neq  \hbar(y_k)$ then there is nothing to prove. If $\hbar(y_j) =  \hbar(y_k)$, without loss of generality assume that $\hbar(y_j) =  \hbar(y_k) = x_1$. Then
\[
g(x_1) = c_1y_1 + \cdots + c_my_m
\]
with $c_j = c_k  = 2$. So
\[
\Sp^1x_1 = x_1^2 + 2y_j + 2y_k \mod \text{other } y\text{-terms}.
\]
Recall that
\[
\Sp^3(y_j) = 2y_j^2x_1 + y_j\lambda(y_j) + \sum_{l, l'}\mu^{(j)}_{l, l'}y_ly_{l'}.
\]
Consequently,
\begin{align*}
y_j^3 = \Sp^1\Sp^3y_j &= \Sp^1\left( 2y_j^2x_1 + y_j\lambda(y_j) + \sum_{l, l'}\mu^{(j)}_{l, l'}y_ly_{l'} \right) \\
&= 2y_j^2(x_1^2+2y_j+2y_k+\cdots) + 2x_1\cdot 2y_j\cdot y_jx_1 + \cdots.
\end{align*}
Note that the term $y_j^2y_k$ is a summand of the right-hand side, but does not appear on the left-hand side. This is a contradiction, and the proof of \Cref{last} is now complete.
\end{proof}

We are now ready to prove the main result of this section.

\begin{proof}[Proof of \Cref{thm:lastthm}]
Suppose that $\chi(\G)\leq n$. Under this assumption, the realizability of $\B(n, \G)$ is a special case of
\Cref{thm:realeyes}. In fact, our proof of \Cref{thm:realeyes} uses \cite[Theorem 1.3]{Tak24}, which shows that a realization $\Xmr_{\B(n, \G)}$ can be taken as a colimit of products of $\BSp(2)$'s and $\BSp(1)$'s. In particular, $\Xmr_{\B(n, \G)}$ receives a map from $\prod_{i=1}^n \BSp(1)$, whose induced map in mod-$3$ cohomology is the desired projection. This proves (b), which then necessarily implies (c). 

\vspace{1mm}

It remains to prove that (c) implies (a). Suppose we have (c), then the standard $\Acal_3$-action on $\ZZ/3[x_1, \cdots, x_n] \cong \Hmr^*(\BSp(1)^{\times n})$ implies \eqref{eq:assume}. So if $\G$ satisfies that $2\G = \G$, then \Cref{last} can be applied to deduce that $\chi(\G)\leq n$, which is (a).

\vspace{1mm}

We now eliminate the assumption that $2\G = \G$ in the implication (c) $\Rightarrow$ (a).
Since the result holds trivially when $2\G = \emptyset$, we assume $2\G \neq \emptyset$.
To simplify notations, we write $\K_{2\G}$ for $\sfdel^{n-1} * 2\G$, and $\K_{\G}$ for $\sfdel^{n-1} * \G$.
Note that if $\sigma$ is a maximal simplex in $\K_{2\G}$, then it is necessarily maximal in $\K_{\G}$.
It follows that there is a map
\[
\lim_{\sigma \in P_{\max}(\K_{\G})} \ZZ/3[\sigma] \longrightarrow \lim_{\sigma \in P_{\max}(\K_{2\G})} \ZZ/3[\sigma]
\]
which is a morphism of $\Acal_3$-algebras. Consider the commutative diagram:
\begin{center}
\begin{tikzcd}
 \B(n, \G) = \lim_{\sigma \in P_{\max}(\K_{\G})} \ZZ/3[\sigma]  \ar[d, "\pi_1"'] \ar[r, "q"] & \B(n, 2\G) = \lim_{\sigma \in P_{\max}(\K_{2\G})} \ZZ/3[\sigma] \ar[d, "\pi_2"] \\
 \ZZ/3[x_1, \cdots, x_n] \ar[r, "="'] & \ZZ/3[x_1, \cdots, x_n]
\end{tikzcd}
.
\end{center}
Suppose that (c) holds for $\G$, so the arrow $\pi_1$ in the above diagram is a map of $\Acal_3$-algebras. We claim that so is the arrow $\pi_2$. Indeed, if $x \in \B(n, 2\G)$, then one can pick a lift $y$ of $x$ in $\B(n, \G)$, as the horizontal arrow $q$ is surjective. Then for any $a\in \ZZ_{>0}$:
\[
\pi_2(\Sp^a x) = \pi_2(\Sp^a q(y)) = \pi_2(q(\Sp^a y)) = \pi_1(\Sp^a y) = \Sp^a\pi_1(y) = \Sp^a\pi_2(q(y)) = \Sp^a\pi_2(x). 
\]
Since the implication (c) $\Rightarrow$ (a) holds for $2\G$, it follows that $\chi(2\G) \leq n$.
If $n\geq 2$, then $\chi(\G) = \chi(2\G)$. Therefore $\chi(\G) \leq n$, and we are done.
Suppose that $n=1$. Because (c) holds for $\G$, \Cref{thm:main3} already implies that $s_3\chi(\G) \leq n = 1$.
Note that when $s_3\chi(\G) = 1$, we also have $\chi(\G) = 1$. This completes the proof.
\end{proof}

\section{The Stanley-Reisner algebras \texorpdfstring{$\A_p(\vb{s}, \G)$}{AsG}} \label{sec:primep}

The algebra $\B(n, \G)$ discussed in \Cref{sec:prime3} admits generalizations at a general prime $p$.
In this section, we introduce, for every odd prime $p$, a family of Stanley-Reisner algebras $\A_p(\vb{s}, \G)$.
These algebras simultaneously generalize $\B(n, \G)$ and the $\A(n, \G)$ studied in \cite{ST25}, and is modeled on classifying spaces $\BSU(k)$ (or a mix of $\BSU(k)$ and $\BSp(l)$).

\vspace{1mm}

Let $p$ be an arbitrarily fixed odd prime. 
We define a family of algebras $\A_p(\vb{s}, \G)$, where 
\begin{itemize}
 \item[-] $\vb{s} = (s_1, \cdots, s_{p-1})$ is a $(p-1)$-tuple of non-negative integers, and
 \item[-] $\G$ is a graph with $m$ vertices,
\end{itemize}
and relate the admission of an $\Acal_p$-structure on $\A_p(\vb{s}, \G)\otimes \ZZ/p$ to the span coloring of $\G$.

\begin{defn} \label[defn]{ApsG}
For each odd prime $p$, we define $\A_p(\vb{s}, \G)$ as the Stanley-Reisner ring $\SR(\K, \upphi)$, where
\[
\K:= {\sf \Delta}^{s_1-1} * \cdots * {\sf \Delta}^{s_{p-1} - 1} * \G,
\]
and $\upphi$ assigns degree $2k+2$ to generators $x_1^{(k)}, \cdots, x_{s_k}^{(k)}$ from ${\sf \Delta}^{s_k-1}$, $k = 1, \cdots, p-1$, and degree $2p+2$ to generators $y_1, \cdots, y_m$ from the graph $\G$. In other words,
\[
\A_p(\vb{s}, \G) = \ZZ[x_1^{(1)}, \cdots, x_{s_1}^{(1)}, \cdots \cdots, x_1^{(p-1)}, \cdots, x_{s_{p-1}}^{(p-1)}, y_1, \cdots, y_m]/I_{\K}.
\]
\end{defn}

Note that $I_{\K}=(y_iy_j)_{(i,j)\not\in \Emr(\G)} \cup (y_iy_jy_k)_{ i<j<k}$. This follows since the minimal missing faces of $M\ast L$ are the unions of the minimal missing faces of $M$ and $L$. Since a simplex has no missing faces, the missing faces of $\K$ above are those of $\G$ and hence we get the description of $I_{\K}$. 

\vspace{1mm}

\subsection{A coloring obstruction to the realizability of \texorpdfstring{$\A_p(\vb{s}, \G)$}{AnG}} \label{sec:col_obs}

In this subsection we prove:

\begin{thm} \label{thm:mainAp}
If $\A_p(\vb{s}, \G)\otimes \ZZ/p$ admits an action of $\Acal_p$, then $s_1 \leq s_p\chi(\G)$.
\end{thm}

We now establish the needed lemmas.

\begin{lemma} \label[lemma]{idealp}
For every $a\geq 0$, $\Sp^a(y_i)$ lies in the ideal
\[
(y_i) + (y_jy_k: 1\leq j < k \leq m) \subset \A_p(\vb{s}, \G).
\]
In particular, $\Sp^p(y_i)$ is of the form
\[
y_i^{p-1}g(y_i) + h(y_i) + \sum_{1\leq j<k \leq m}\mu^{(i)}_{j, k}y_jy_k,
\]
where 
\begin{itemize}
\item[-] $g(y_i)$ is a $\ZZ/p$-linear combination of $x_1^{(1)}, \cdots, x_{s_1}^{(1)}$.
\item[-] $h(y_i)$ is an element in $(y_i)$ whose summands are divisible at most by $y_i^{p-2}$.
\item[-] $\mu^{(i)}_{j, k}$ is some $\ZZ/p$-linear combination of monomials in $x_1^{(1)}, \cdots, x_{s_{p-1}}^{(p-1)}$.
\end{itemize}
\end{lemma}

\begin{proof}
The first part can be proved in a similar way as in \Cref{square_formulas}, using \Cref{preserve} and \Cref{ideal}.
For the expression of $\Sp^p(y_i) $, note that it belongs to $(y_i) + (y_jy_k: 1\leq j < k \leq m)$, and that $|y_i| = 2p+2$ and $\Sp^p$ raises degree by $2p(p-1)$. Consequently, $\Sp^p(y_i)$ has degree $2p^2+2$, which forces $\Sp^p(y_i) $ to be of the desired form.
\end{proof}

The values $g(y_i)$ can be utilized to define a function
\[
g: \{y_1, \cdots, y_m\} \longrightarrow \ZZ/p \, \langle x_1^{(1)}, \cdots, x_{s_1}^{(1)} \rangle.
\]
We claim that $g$ gives rise to an $s_1$-coloring.

\begin{lemma} \label[lemma]{A_2G_coloring}
Suppose that every vertex of $\G$ has degree at least $2$, and that $\A_p(\vb{s}, \G)\otimes \ZZ_p$ admits an action of $\Acal_p$. Then for every $i$, the linear inclusion
\[
\langle g(y_j) : y_j\in N(y_i) \rangle \longrightarrow \langle g(y_j) : y_j\in N(y_i) \cup \{y_i\} \rangle
\]
has nonzero cokernel.
\end{lemma}

\begin{proof}
Assume for contradiction that the cokernel is zero. Then $g(y_i)$ is a linear combination 
\begin{equation} \label{eq:assumptionp}
g(y_i) = \sum_{y_j\in N(y_i)} a_jg(y_j), \,\, a_j\in \ZZ/p.
\end{equation}
We now apply the Adem relation $\Sp^1\Sp^p = \Sp^{p+1}$ to $y_i$, and use \Cref{idealp}. One obtains
\begin{equation} \label{eq:AdemPp}
y_i^p = \Sp^{p+1}(y_i) = \Sp^1(\Sp^p(y_i)) = \Sp^1\left(y_i^{p-1}g(y_i) + h(y_i) + \sum_{1\leq j < k \leq m}\mu^{(i)}_{j, k}y_jy_k\right).
\end{equation}
The right-hand side of \Cref{eq:AdemPp} must contain $y_i^p$ as a summand, dictating that $\Sp^1(g(y_i))$ contains $y_i$ as a summand. Now applying $\Sp^1$ to \Cref{eq:assumptionp}, one concludes that there exists some $j$ with $y_j\in N(y_i)$, so that $\Sp^1(g(y_j))$ contains $y_i$ as a summand as well. For such $j$, we find that
\begin{equation} \label{eq:AdemPpj}
y_j^p = \Sp^{p+1}(y_j) = \Sp^1\Sp^p(y_j) = \Sp^1\left( y_j^{p-1}g(y_j) + h(y_j) + \sum_{1\leq k < l \leq m}\mu^{(j)}_{k, l}y_ky_l \right)
\end{equation}
Since $\Sp^1(g(y_j))$ contains $y_i$ as a summand, the right-hand side of \Cref{eq:AdemPpj} contains $y_j^{p-1}y_i$ as a summand. Note that $y_j^{p-1}y_i$ is nonzero as $y_j\in N(y_i)$. However, $y_j^{p-1}y_i$ does not appear in the left-hand side of \Cref{eq:AdemPpj}, which is a contradiction.
\end{proof}

It follows from \Cref{A_2G_coloring} that $s_p\chi(2\G) \leq n$, using \cite[Lemma 6.10 and Lemma 6.11]{ST25}. We now prove the main result, \Cref{thm:mainAp}.

\begin{proof}[Proof of \Cref{thm:mainAp}]
We discuss two cases. Suppose that $2\G \neq \emptyset$. Since $s_p\chi(2\G) \leq n$, there is map $\lambda: 2\G\longrightarrow \sfA_{(\ZZ/p)^n}$ from $2\G$ to the graph $ \sfA_{(\ZZ/p)^n}$ representing $n$-span colorings over $\ZZ/p$. 
Note that there is a finite sequence of graph inclusions $2\G=\G_N\subset \G_{N-1} \subset \cdots \subset \G_1 \subset \G_0=\G$, where each step adds a vertex of degree 0 or 1. 
By \cite[Theorem 6.12]{ST25}, any map from $\G_{i+1}$ to $\sfA_{(\ZZ/p)^n}$ can be extended over $\G_i$. So $\lambda: 2\G\longrightarrow \sfA_{(\ZZ/p)^n}$ can eventually be extended to a map $\tilde{\lambda}: \G \longrightarrow \sfA_{(\ZZ/p)^n}$. This corresponds to an $n$-span coloring of $\G$, and consequently $s_p\chi(\G) \leq n$.

\vspace{1mm}

Now suppose that $2\G = \emptyset$. In this case $\G$ is a union of trees, and in particular $s_p\chi(\G) \leq 2$. So if $s_1\geq 2$ then automatically $s_p\chi(\G) \leq 2 \leq s_1$. It then suffices to prove the result for $s_1 = 0$ and $s_1 = 1$.

\vspace{1mm}

When $s_1=0$, we need to show that $s_p\chi(\G) = 0$, i.e., $\G$ is empty. Assume for contradiction that $\G \neq \emptyset$. Let $\sigma$ be a maximal simplex in $\K$. By \cite[Proposition 6.4]{ST25}, the algebra $\ZZ/p[\sigma]$ has an induced $\Acal_p$-action. Note that $\ZZ/p[\sigma]$ is either the algebra $A$ or $A'$ below:
\begin{itemize}
 \item $A:= \ZZ/p[x_1, \cdots, x_a, y]$, with $|x_i| \in \{6, 8,  \cdots, 2p\}$ for $1\leq i \leq a$, and $|y| = 2p+2$;
 \vspace{1mm}
 \item $A':= \ZZ/p[x_1, \cdots, x_a, y_1, y_2]$, with $|x_i| \in \{6, 8,  \cdots, 2p\}$ for $1\leq i \leq a$, and $|y| = 2p+2$.
\end{itemize}
We show that $A$ or $A'$ cannot possibly support an $\Acal_3$-action. For the algebra $A$, consider applying the Adem relation $\Sp^1\Sp^p = \Sp^{p+1}$ to $y$. 
Since none of the $x_i$'s is of degree 4, the monomials appearing in the expression of $\Sp^p(y)$ (which are of degree $2p^2+2$) cannot be divisible by $y^{p-1}$ (which has degree $2p^2-2$), or any higher power of $y$. 
So $\Sp^p(y)$ is of the form
\[
y^{p-2}\lambda(y) + \mu(y),
\]
where terms in $\lambda(y)\in A$ cannot be divisible by $y$, and terms in $\mu(y)\in A$ can at best be divisible by $y^{p-3}$. However, the identity
\[
y^p = \Sp^{p+1}(y) = \Sp^1\Sp^p(y) = \Sp^1 \left( y^{p-2}\lambda(y) + \mu(y) \right)
\]
gives rise to a contradiction, since the right-hand side cannot contain $y^p$ as a summand. A similar argument (applying the same Adem relation to $y_1$ or $y_2$) shows that $A'$ cannot support an $\Acal_p$-action.

\vspace{1mm}

When $n=1$, we need to show that $s_p\chi(\G) \leq 1$, i.e., $\G$ contains no edges.
Assume not. Let $\sigma$ be a maximal simplex in $\K$. Again by \cite[Proposition 6.4]{ST25}, there is an induced $\Acal_p$-action on $\ZZ/3[\sigma]$, which is the algebra
\begin{itemize}
 \item $B:=\ZZ/p[x_1, \cdots, x_b, y_1, y_2]$, with $|x_1| = 4$, $|x_i| \in \{6, 8,  \cdots, 2p\}$ for $2\leq i \leq b$, and $|y_1| = |y_2| = 2p+2$.
\end{itemize}
Now that there is a single degree $4$ generator $x_1\in B$, the element $\Sp^p(y_1)$ can be written as
\[
c_1x_1y_1^{p-1} + \xi_1(y_1),
\]
where $c_1\in \ZZ/p$ is a constant, and $\xi_1(y_1)\in B$ is a polynomial whose $y_1$-degree is no greater than $p-2$. Similarly, there is an expression
\[
\Sp^p(y_2) = c_2x_1y_2^{p-1} + \xi_2(y_2),
\]
where $c_2\in \ZZ/p$ is a constant, and $\xi_2(y_2)\in B$ is a polynomial whose $y_2$-degree is no greater than $p-2$.
Applying the Adem relation $\Sp^1\Sp^p = \Sp^{p+1}$ to $y_1$, one obtains that
\begin{align*}
y_1^p = \Sp^{p+1}y_1 = \Sp^1\Sp^p y_1 &= \Sp^1\left(c_1x_1y_1^{p-1} + \xi_1(y_1)\right) \\
 &= c_1\Sp^1(x_1)y_1^{p-1} + c_1(p-1)x_1y_1^{p-2}\Sp^1(y_1) + \Sp^1 \xi_1(y_1).
\end{align*}
The above equality forces $c_1 = 1$, and $\Sp^1(x_1) = y_1$. However, the same calculation with $y_2$ would imply that $c_2 = 1$ and $\Sp^1(x_1) = y_2$, yielding a contradiction.
\end{proof}

\vspace{1mm}

As an immediate consequence of  \Cref{thm:mainAp}, we obtain the following:

\begin{cor} \label[cor]{ARealize}
If $\A_p(\vb{s}, \G)$ is realizable, then $s_p\chi(\G) \leq s_1$.
\end{cor}

\vspace{1mm}



\subsection{The realizability of \texorpdfstring{$\A_p(\vb{s}, \G)$}{AnG}} \label{sec:top_inv}

In this subsection, we present some sufficient conditions for the realizability of $\A_p(\vb{s}, \G)$. We will 
assume that $\vb{s}=\vb{r}+\vb{t}$ where 
$\vb{r}=(r_1,\cdots,r_{p-1})$ and $\vb{t}=(t_1,\cdots,t_{p-1})$ are sequences such that $r_i\geq r_{i+1}$, 
$t_{2i}=0$ and $t_i\geq t_{i+2}$. The first sequence corresponds to products of various $\BSU(n)$ and the second to products of $\BSp(n)$. 
The goal is to prove that if $r_{p-1}+t_{p-2}\geq \chi(\G)$ then $\A_p(\vb{r}+\vb{t}, \G)$ is realizable.
In fact, we can prove a bit more, namely the realizability of $\A_p(\vb{r}+\vb{t}, L)$ where $L$ is a simplicial complex with $L^1 = \G$ (see \Cref{thm:realeyes}).

\begin{defn} \label[defn]{ApsL}
Let $L$ be a simplicial complex with vertex set $\V(L) = \{y_1, \cdots, y_m\}$. We define $\A_p(\vb{s}, L)$ as the Stanley-Reisner algebra $\SR(K, \upphi)$, where
\[
\K= {\sf \Delta}^{s_1-1} * \cdots * {\sf \Delta}^{s_{p-1} - 1} * L,
\]
and $\upphi$ assigns degree $2k+2$ to generators $x_1^{(k)}, \cdots, x_{s_k}^{(k)}$ from ${\sf \Delta}^{s_k-1}$, $k = 1, \cdots, p-1$, and degree $2p+2$ to generators $y_1, \cdots, y_m$ from the vertices $L$. 
\end{defn}

We use the following result of Takeda \cite{Tak24}:

\begin{thm}[\cite{Tak24}, Theorem 1.3] \label{thm:takeda1.3} 
Let $\SR(\K,\phi)$ be the graded Stanley-Reisner ring for simplicial complex $\K$ with the vertex set $\V$ and $\phi : \V \rightarrow 2\ZZ_{>0}$. We can construct a space $X$ as a homotopy colimit such that $\Hmr^*(X) = \SR(\K, \phi)$, if $\SR(\K, \phi)$ satisfies the following condition:

\begin{itemize}
\item There is a partition of $\V$ , $\coprod A_i=\V$, such that for all $i$ and $\sigma\in \Pmax(\K)$
the set $\{\phi(x) | x\in \sigma\cap A_i\}$ is equal to $\{2, 2, . . . , 2\}$, $\{4, 6, . . . 2n\} \cup \{2, 2, . . . , 2\}$
or $\{4, 8, . . . 4n\} \cup \{2, 2, . . . , 2\}$ as a multiset for some $n$.
\end{itemize}
\end{thm}

\begin{thm} \label{thm:realeyes}
Let $L$ be a simplicial complex with $L^1 = \G$. If $r_{p-1}+t_{p-2}\geq \chi(\G)$ then $\A_p(\vb{r}+\vb{t}, L)$ is realizable. 
\end{thm}

\begin{proof}
We wish to partition all the generators of $\A_p(\vb{s}, L)$ - first recalling (see \Cref{ApsL}) that the $s_i$ gives the number of generators $x$ with 
$|x|=2(i+1)$ and the generators $y$ with degree $2p+2$ correspond to the vertices of the simplicial complex $L$.

\vspace{1mm}
For $1\leq i\leq p-2$ and $j\in r_i-r_{i+1}$ let $U(i,j)$ be a set of generators with multidegrees
$[4, 6, \cdots, 2i+2]$ and for $1\leq i \leq \frac{p-3}{2}$ take $W(i,k), k\in [t_{2i-1}-t_{2i+1}]$ a set of generators with multidegrees $[4,8, \cdots, 4i]$. We also ask that the $U(i,j)$ and $W(i,k)$ are all pairwise disjoint. 

\vspace{1mm}
The definitions of $U(p-1, j)$  and $W((p-1)/2, k)$ for $j\in [r_{p-1}]$ and $k\in [t_{p-1}]$
are slightly different. Since $\chi(L^1)\leq r_{p-1}+t_{p-1}$ there is a map
$\phi$ from  $V(L)$ to the set of the $U(p-1, j)$  and the $W((p-1)/2, k)$ so that no edge of $L^1$ 
has vertices labeled by the same set. The generators inside $U(p-1, j)$ would have the multidegrees 
$[4, 6 \cdots,  2p]\bigcup \cup_{\phi^{-1}(U(p-1, j)} [2p+2]$. The first generators being the ones 
as in the other cases and the second bunch of generators corresponding to the vertices of $L^1$ whose
colour is $U(p-1,j)$ (or $W((p-1)/2,k)$). 

\vspace{1mm}
Clearly the $U(i,j)$ and $W(i,k)$ give a partition of the generators. Next let  
$\sigma\in \Pmax(L)$, then each generator of degree $<2p+2$ (so the $x$s) is in $\sigma$. 
If follows that if $i<p-1$ then $\sigma\cap U(i,j)=U(i,j)$. 
In addition $\sigma$ has some number of degree $2p+2$ generators (so the $y$s). Intersecting $\sigma$ 
with any colour $\sigma\cap U(p-1,j)$ (or $\sigma \cap W((p-1)/2,j)$ there can be at most one 
$y_l$ since if $y_l,y_{l'}\in U(n-1,j)$ there is no edge between them and so they can't be together in any simplex. 
Thus $\phi(\sigma\cap U(p-1,j))=\{4, \cdots, 2p\}$ or $\{4, \cdots 2p, 2p+2\}$ and 
$\sigma\cap W((p-1)/2,j)=\{4, \cdots, 2(p-1)\}$ or $\{4, \cdots 2(p-1), 2p+2\}$.

\vspace{1mm}
We have just shown the hypotheses of \Cref{thm:takeda1.3} and so we can conclude that $\A_p(\vb{s}, L)$ is realizable. 
\end{proof}

\section{Interpreting realization as a graph invariant}\label{OOG}

In this section we look at special constant sequences $\vb{r}(n)$ and $\vb{t}(n)$. They are defined by $r(n)_i=n$ for 
$1\leq i\leq p-1$ and 
$t(n)_{2i-1}=n$ when $1\leq i \leq (p-1)/2$ (as before $t(n)_{2i}=0$). 
Note that $A_2(\vb{r}(n), \G)=A(n,G)$ and 
$A_3(\vb{t}(n), \G)=B(n,G)$. 
For these sequences we have both the upper 
bound of Theorem \ref{thm:realeyes} and the lower bound of Theorem \ref{thm:mainAp}.  
We introduce topologically defined graph invariants so we can interpret 
these theorems as sequences of inequalities of graphs invariants which are analogues of  Theorem \ref{ST3} for $p>2$. 

\begin{defn} \label[defn]{topchrom}
Let $\G$ be a graph. 
The $\A$-topological chromatic number of $\G$, which we write as $\chi_{\Top_p, \A}(\G)$, is the smallest $n$ 
so that $\A_p(\vb{r}(n), \G)$ is realizable. The $\B$-topological chromatic number of a 
graph $\G$, denoted by $\chi_{\Top_p, \B}(\G)$, is the smallest $n$ so that $\B_p(\vb{t}(n), \G)$ is realizable. 

\end{defn}

There are also algebraically defined ones:

\begin{defn}
For a graph $G$ and a prime $p$. Define $\chi_{Steen_p, A}$ as the least $n$ such that 
$\A_p(\vb{r}(n), \G)\otimes \ZZ/p$ has an $\mathcal A_p$ action and if $p\geq 3$, define 
$\chi_{Steen_p, B}$ as the least $n$ such that 
$\A_p(\vb{t}(n), \G)\otimes \ZZ/p$ has an $\mathcal A_p$ action.

\end{defn}

\begin{thm} \label{thm:chronumber} For $p>2$ 
the following comparisons of chromatic numbers hold:
\[
s_p\chi(\G) \leq \chi_{Steen_p, A}\leq \chi_{\Top_p, \A}(\G) \leq \chi(\G), \quad \text{and} 
\quad s_p\chi(\G) \leq \chi_{Steen_p, B}\leq \chi_{\Top_p, \B}(\G) \leq \chi(\G).
\]

\end{thm}

\begin{proof}
The inequalities $s_p\chi(\G) \leq \chi_{\Top, \A}(\G)$ and $s_p\chi(\G) \leq \chi_{\Top, \B}(\G)$ follow from 
Theorem \ref{thm:mainAp}. The middle inequalities follow since if any $A\cong H^*(X, \ZZ)$ is realizable then 
$A\otimes \ZZ/p$ has an $\mathcal A_p$ action. Finally Theorem \ref{thm:realeyes} implies $\chi_{\Top_p, \A}(\G) \leq \chi(\G)$ and 
$\chi_{\Top_p, \B}(\G) \leq \chi(\G)$.
\end{proof}

It was shown in \cite{ST25} that for every $p$ there is a $\G$ such that $s_p\chi(\G)<\chi(\G)$. We also know that 
$s_p\chi (A_{(Z/p)^n})=n$ and it was recently shown \cite{ILMS} that 
\[
\lim\limits_{n\rightarrow \infty} \frac{ \log \chi (A_{(\ZZ/p)^{\times n}}) }{n} = \log (p),
\] 
so the two middle invariants have a large gap to sit in. Theorem \ref{ST1} says 
$\chi_{Steen_2, A}=\chi_{\Top_2, \A}(\G)$, but we do not know more about which of these inequalities 
can be strict or where in the large gap these invariants sit. Theorem \ref{thm:lastthm} above might indicate for 
$p=3$ the invariants are closer to 
$\chi(G)$.


\bibliographystyle{amsalpha}

\bibliography{SROP}

\end{document}